\documentclass{article}

\usepackage{amsmath,amssymb,amsthm}

\newtheorem{theorem}{Theorem}[section]

\theoremstyle{remark}
\newtheorem{remark}{Remark}[section]
\theoremstyle{definition}
\newtheorem{definition}{Definition}[section]

\theoremstyle{definition}

\begin{document}

\title{A Noether Theorem on Unimprovable Conservation Laws
for Vector-Valued Optimization Problems in Control
Theory\footnote{To appear (paper ID: 131) in the Conference Proceedings
of ICMSAO'05 -- First International Conference on Modeling,
Simulation and Applied Optimization, February 1-3 2005,
American University of Sharjah, United Arab Emirates.}}

\author{Delfim F. M. Torres\\
        \texttt{delfim@mat.ua.pt}}

\date{Control Theory Group (\textsf{cotg})\\
      Centre for Research in Optimization and Control\\
      Department of Mathematics, University of Aveiro\\
      3810-193 Aveiro, Portugal}

\maketitle


\begin{center}
Dedicated to the memory of Almaskhan Gugushvili
\end{center}

\bigskip


\begin{abstract}
We obtain a version of Noether's invariance theorem
for optimal control problems with a finite number of cost functionals.
The result is obtained by formulating E.~Noether's result to optimal
control problems subject to isoperimetric constraints, and then
using the unimprovable (Pareto) notion of optimality. A result of this kind
was posed to the author, as a mathematical open question,
of great interest in applications of control engineering, by A.~Gugushvili.
\end{abstract}


\vspace*{0.5cm}



\vspace*{0.3cm}



\section{Introduction}

E. Noether's Theorem, which relates symmetries and conservation laws,
is one of the most deep and rich, powerful and helpful results of physics.
It describes the fundamental fact that ``invariance with respect
to some group of parameter transformations gives rise
to the existence of conservation laws''.
Typical application of conservation laws is
to lower the order of the systems.
They are, however, a useful tool for
many other reasons, \textrm{e.g.},
to prove regularity of the minimizers
in the calculus of variations and optimal control \cite{kiev2003}.
Noether's theorem comprises all results
on conservation laws known to classical mechanics.
Thus, for instance, the invariance relative to translation
with respect to time yields conservation of energy; while
conservation of linear and angular momenta reflect, respectively,
translational and rotational invariance.
Noether's theorem is applicable also in quantum mechanics,
field theory, electromagnetic theory,
and has deep implications in the general theory of relativity.
It is useful to explain a myriad of things,
from the fusion of hydrogen to the motion of planets
orbiting the sun \cite{MR53:10537}. Moreover, it turns out that
Noether's theorem is much more than a theorem:
it is a \emph{principle}, which can be formulated, as a theorem,
in many different contexts, under many different assumptions.
It is possible, \textrm{e.g.}, to formulate the classical Noether's theorem
of the calculus of variations for bigger classes of nonsmooth
admissible functions \cite{cpaa}; in the more general
context of optimal control \cite{alik,ejc}; or to obtain
discrete-time versions \cite{discreteCAO03}.
For an account of Noether's symmetry principle
in the context of optimal control,
the use of conservation laws
to integrate and decrease the order
of the equations given by the Pontryagin
maximum principle \cite{MR29:3316b},
and for practical examples,
such as the problem of synchronization
of difficult control systems,
we refer the reader to \cite{alik}.
Here we are interested in generalizing the previous results
to cover optimal control problems which, in place
of a single cost functional, have a vector-valued functional to minimize.
For an introduction to problems of optimal control
with multiple objectives, we refer the reader to Salukvadze's book
\cite{MR81k:93001}. Multiobjective optimal control attracts more and more attention,
and is source of many open questions \cite{AndreasControlo04}.
The motivation for the present study was a
challenge proposed to the author by A.~Gugushvili
on November 18, 2003. A.~Gugushvili wanted to generalize the symmetry and
conservation laws to multiobjective problems of optimal control:
\emph{``We would like to develop E.~Noether's theory for
multicriteria optimal control systems. If you have
any ideas and work on these problems, please, let us know.''}
Theorem~\ref{r:MC:consLaw:PO} is, to the best of our
knowledge, the first attempt in this direction.


\section{Optimal Control with Isoperimetric Constraints}

It is well known that necessary optimality conditions for optimal control
problems subject to isoperimetric constraints, are also necessary
for unimprovable (Pareto) optimality in the problem with a vector-valued cost
(\textrm{cf.}, \textrm{e.g.},  \cite[Chap.~17]{MR84m:49002},
\cite[p.~22]{MR81k:93001}). Consider a nonlinear control system,
\begin{equation}
\label{eq:ncs}
\dot{x}(t) = \varphi\left(t,x(t),u(t)\right) \, ,
\end{equation}
of $n$ differential equations, subject to $k$ isoperimetric
equality constraints,
\begin{equation}
\label{eq:ec:ig}
\int_a^b g_i\left(t,x(t),u(t)\right) \mathrm{d}t = \xi_i \, , \quad i = 1,\ldots, k \, ;
\end{equation}
$m$ isoperimetric inequality constraints,
\begin{equation}
\label{eq:ec:des}
\int_a^b g_j\left(t,x(t),u(t)\right) \le \xi_j \, , \quad j = k+1,\ldots, k+m \, ;
\end{equation}
and $2n$ boundary conditions
\begin{equation}
\label{eq:bc}
x(a) = \alpha \, , \quad x(b) = \beta \, .
\end{equation}
The problem consists to find a piecewise-continuous control function
$u(\cdot) = \left(u_1(\cdot),\ldots,u_r(\cdot)\right)$, taking
value on a given set $\Omega \subseteq \mathbb{R}^r$, and the
corresponding state trajectory $x(\cdot) =
\left(x_1(\cdot),\ldots,x_n(\cdot)\right)$, satisfying
\eqref{eq:ncs}, \eqref{eq:ec:ig}, \eqref{eq:ec:des}, and
\eqref{eq:bc}, which minimizes (or maximizes) the (scalar)
integral cost functional
\begin{equation*}
I[x(\cdot),u(\cdot)] = \int_a^b L(t,x(t),u(t)) \mathrm{d}t \, .
\end{equation*}
This problem is denoted in the sequel by $(P_1)$.
Both the initial time $a$ and terminal time $b$, $a < b$, are fixed.
The boundary values $\alpha,\,\beta \in \mathbb{R}^n$,
and constants $\xi_i$, $i = 1,\ldots,k+m$, are also
given. The functions $L(\cdot,\cdot,\cdot)$, $\varphi(\cdot,\cdot,\cdot)$
and $g(\cdot,\cdot,\cdot)$ are assumed to be continuously differentiable
with respect to all variables.
The celebrated Pontryagin's maximum principle \cite{MR29:3316b}
gives necessary optimality conditions to be satisfied
by the solutions of optimal control problems.
Formulation of the maximum principle for
problems with isoperimetric constraints can be found,
\textrm{e.g.}, in \cite[\S 13.12]{MR84m:49002}.
\begin{theorem}[Pontryagin Maximum Principle for $(P_1)$]
\label{PMP:P1}
Let $u(t)$, $t \in [a,b]$, be an optimal control for
the isoperimetric (scalar) optimal control problem $(P_1)$, and
$x(\cdot)$ the corresponding state trajectory. Then there
exists a constant $\psi_0 \le 0$, a continuous costate
$n$-vector function $\psi(\cdot)$ having piecewise-continuous derivatives,
and constant multipliers
$\lambda_i$, $i = 1, \ldots, k + m$,
where $(\psi_0,\psi(\cdot),\lambda) \neq 0$,
satisfying the \emph{pseudo-Hamiltonian system}
\begin{equation*}
\begin{cases}
\dot{x}(t) &=
\frac{\partial H}{\partial \psi}\left(t,x(t),u(t),\psi_0,\psi(t),\lambda\right)
\, , \\
\dot{\psi}(t) &= -
\frac{\partial H}{\partial x}\left(t,x(t),u(t),\psi_0,\psi(t),\lambda\right)
\, ;
\end{cases}
\end{equation*}
the \emph{maximality condition}
\begin{equation*}
H\left(t,x(t),u(t),\psi_0,\psi(t),\lambda\right) = \max_{u \in \Omega}
H\left(t,x(t),u,\psi_0,\psi(t),\lambda\right) \, ;
\end{equation*}
where the Hamiltonian $H$ is defined by
\begin{equation}
\label{eq:Hamiltonian:P1}
H(t,x,u,\psi_0,\psi,\lambda)
= \psi_0 L(t,x,u) + \psi \cdot \varphi(t,x,u) + \lambda \cdot g(t,x,u)\, .
\end{equation}
Moreover, $\lambda_j \le 0$, $j = k+1,\ldots, k+m$, where $\lambda_j = 0$ if
\begin{equation*}
\int_a^b g_j\left(t,x(t),u(t)\right) < \xi_j \, ;
\end{equation*}
and $H(t,x(t),u(t),\psi_0,\psi(t),\lambda)$ is a continuous
function of $t$ and, on each interval of continuity of $u(\cdot)$,
is differentiable and satisfies the equality
\begin{equation}
\label{eq:dHdt}
\frac{dH}{dt}\left(t,x(t),u(t),\psi_0,\psi(t),\lambda\right)
= \frac{\partial H}{\partial t}\left(t,x(t),u(t),\psi_0,\psi(t),\lambda\right) \, .
\end{equation}
\end{theorem}


\section{Vector-Valued Optimal Control Problems}
\label{sec:VecValOCP}

When optimal control is used to model a real problem,
it is natural that several (conflicting) cost functionals (``objectives'')
are desired to be taken in account (see \cite{MR81k:93001}
for many practical situations). The problem is then
to minimize a vector-valued functional with components
\begin{equation*}
I_i[x(\cdot),u(\cdot)] = \int_a^b L_i(t,x(t),u(t)) \mathrm{d}t \, ,
\quad i = 1,\ldots,N \, ,
\end{equation*}
subject to a dynamical control system \eqref{eq:ncs},
and boundary conditions \eqref{eq:bc}. We denote this
problem by $(P)$.
\begin{definition}
\label{def:PO}
An admissible pair $\left(\tilde{x}(\cdot),\tilde{u}(\cdot)\right)$
is said to be an \emph{unimprovable solution}, \emph{compromise solution},
or a \emph{Pareto solution} for $(P)$ if, and only if, for every admissible
pair $\left(x(\cdot),u(\cdot)\right)$, either
\begin{equation*}
I_i[\tilde{x}(\cdot),\tilde{u}(\cdot)] =
I_i[x(\cdot),u(\cdot)] \quad
\forall i \in \{1,\ldots,N\} \, ,
\end{equation*}
or there exists at least one $i \in \{1,\ldots,N\}$ such that
\begin{equation*}
I_i[\tilde{x}(\cdot),\tilde{u}(\cdot)] < I_i[x(\cdot),u(\cdot)] \, .
\end{equation*}
\end{definition}
It turns out that necessary conditions
for optimal control problems with isoperimetric constraints,
are also necessary for Pareto-optimality of optimal control
problems with a vector-valued cost. Theorem~\ref{th:NCINCPO}
is a simple consequence of Definition~\ref{def:PO}
(\textrm{cf.}, \textrm{e.g.}, \cite[Theorem~17.1]{MR84m:49002}).
\begin{theorem}
\label{th:NCINCPO}
If $\left(\tilde{x}(\cdot),\tilde{u}(\cdot)\right)$ is a Pareto-solution
of problem $(P)$, then it is a minimizer for the isoperimetric
optimal control problems with the integral scalar-valued cost
\begin{equation*}
I_i[x(\cdot),u(\cdot)] \, , \quad i \in \{1,\ldots,N\} \, ,
\end{equation*}
and isoperimetric constraints
\begin{equation*}
I_j[x(\cdot),u(\cdot)] \le I_j[\tilde{x}(\cdot),\tilde{u}(\cdot)]
\, , \quad j = 1,\ldots,N \text{ and } j \ne i \, .
\end{equation*}
\end{theorem}
From Theorems~\ref{th:NCINCPO} and \ref{PMP:P1}
(Pontryagin maximum principle for problems with
isoperimetric constraints) it follows
the so called ``general theorem of optimal control''
(\textrm{cf.} \cite[p.~22]{MR81k:93001}).
\begin{theorem}
\label{th:GenThOptCont}
If $\left(\tilde{x}(\cdot),\tilde{u}(\cdot)\right)$ is a Pareto-solution
of problem $(P)$, then there exists a continuous costate
$n$-vector function $\psi(\cdot)$ having piecewise-continuous derivatives,
and constant multipliers $\lambda = \left(\lambda_1,\ldots,\lambda_N\right)$,
where $(\psi(\cdot),\lambda) \neq 0$,
satisfying the \emph{pseudo-Hamiltonian system}
\begin{equation*}
\begin{cases}
\dot{x}(t) &=
\frac{\partial \mathcal{H}}{\partial \psi}\left(t,x(t),u(t),\psi(t),\lambda\right)
\, , \\
\dot{\psi}(t) &= -
\frac{\partial \mathcal{H}}{\partial x}\left(t,x(t),u(t),\psi(t),\lambda\right)
\, ;
\end{cases}
\end{equation*}
the \emph{maximality condition}
\begin{equation*}
\mathcal{H}\left(t,x(t),u(t),\psi(t),\lambda\right) = \max_{u \in \Omega}
\mathcal{H}\left(t,x(t),u,\psi(t),\lambda\right) \, ;
\end{equation*}
where the Hamiltonian $\mathcal{H}$ is defined by
\begin{equation}
\label{eq:Hamiltonian:P}
\mathcal{H}(t,x,u,\psi,\lambda)
= \lambda \cdot L(t,x,u) + \psi \cdot \varphi(t,x,u) \, .
\end{equation}
Moreover, $\lambda_j \le 0$, $j = 1,\ldots, N$;
and $\mathcal{H}(t,x(t),u(t),\psi(t),\lambda)$ is a continuous
function of $t$ and, on each interval of continuity of $u(\cdot)$,
is differentiable and satisfies the equality
\begin{equation*}
\label{eq:dHdt:PO}
\frac{d\mathcal{H}}{dt}\left(t,x(t),u(t),\psi(t),\lambda\right)
= \frac{\partial \mathcal{H}}{\partial t}\left(t,x(t),u(t),\psi(t),\lambda\right) \, .
\end{equation*}
\end{theorem}


\section{Main Results: Noether-type Theorems}

Theorem~\ref{r:MC:consLaw:IC} asserts that
the presence of a symmetry for the optimal control
problems involving equality and inequality isoperimetric constraints,
imply that their Pontryagin extremals (and solutions) preserve
a well-defined quantity (there exists a conservation law associated
with each symmetry). The result is formulated, as it happens for the problems
of the calculus of variations \cite{cpaa},
and for the unconstrained scalar-valued continuous \cite{ejc}
and discrete-time \cite{discreteCAO03} optimal control problems,
as an instance of Noether's universal principle.

\begin{definition}
An equation $C\left(t,x(t),u(t),\psi_0,\psi(t),\lambda\right) = \text{constant}$,
valid in $t \in [a,b]$ for any quintuple
$\left(x(\cdot),u(\cdot),\psi_0,\psi(\cdot),\lambda\right)$
satisfying the Pontryagin maximum principle (Theorem~\ref{PMP:P1}),
is called a \emph{conservation law} for problem $(P_1)$.
\end{definition}

\begin{theorem}[Noether theorem for optimal control problems
with isoperimetric constraints]
\label{r:MC:consLaw:IC}
If there exists a $C^2$-smooth one-parameter group of transformations
\begin{gather*}
h^s : [a,b] \times \mathbb{R}^n \times \mathbb{R}^r \rightarrow
       \mathbb{R} \times \mathbb{R}^n \times \mathbb{R}^r \, , \\
h^s(t,x,u) = \left(T(t,x,u,s), X(t,x,u,s), U(t,x,u,s)\right) \, , \\
s \in (-\varepsilon, \varepsilon)  \, , \, \varepsilon > 0 \, ,
\end{gather*}
with $h^0(t,x,u) = (t,x,u)$ for all
$(t,x,u) \in [a,b] \times \mathbb{R}^n \times \mathbb{R}^r$, and
satisfying
\begin{gather}
L\left(t,x(t),u(t)\right)
= L \circ h^s\left(t,x(t),u(t)\right)
\frac{d}{dt} T\left(t,x(t),u(t),s\right) \, , \label{eq:invi} \\
\frac{d}{dt} X\left(t,x(t),u(t),s\right)
= \varphi \circ h^s\left(t,x(t),u(t)\right)
\frac{d}{dt} T\left(t,x(t),u(t),s\right) \, ,  \label{eq:invii} \\
g\left(t,x(t),u(t)\right)
= g \circ h^s\left(t,x(t),u(t)\right)
\frac{d}{dt} T\left(t,x(t),u(t),s\right)  \, , \label{eq:inviii}
\end{gather}
then,
\begin{multline*}
\psi(t) \cdot \frac{\partial}{\partial s}
\left.X\left(t,x(t),u(t),s\right)\right|_{s = 0} \\
- H(t,x(t),u(t),\psi_0,\psi(t),\lambda) \frac{\partial}{\partial s}
\left.T\left(t,x(t),u(t),s\right)\right|_{s = 0} = \text{const}
\end{multline*}
is a conservation law for problem $(P_1)$,
with $H$ the Hamiltonian \eqref{eq:Hamiltonian:P1}
associated to the problem $(P_1)$.
\end{theorem}

\begin{proof}
Using the fact that $h^0(t,x,u) = (t,x,u)$,
from condition \eqref{eq:invi} one gets
\begin{align}
0 & = \frac{d}{ds}
\left.\left(L \circ h^s\left(t,x(t),u(t)\right)
\frac{d}{dt} T\left(t,x(t),u(t),s\right)\right)\right|_{s = 0} \notag \\
  & = \frac{\partial L}{\partial t}
      \left.\frac{\partial T}{\partial s}\right|_{s = 0}
      + \frac{\partial L}{\partial x} \cdot
      \left.\frac{\partial X}{\partial s}\right|_{s = 0}
      + \frac{\partial L}{\partial u} \cdot
      \left.\frac{\partial U}{\partial s}\right|_{s = 0}
      + L \frac{d}{dt} \left.\frac{\partial T}{\partial s}\right|_{s =
      0} \, , \label{eq:ds0i}
\end{align}
while condition \eqref{eq:invii} and \eqref{eq:inviii} yields
\begin{gather}
\frac{d}{dt} \left.\frac{\partial X}{\partial s}\right|_{s = 0}
= \frac{\partial \varphi}{\partial t}
      \left.\frac{\partial T}{\partial s}\right|_{s = 0}
      + \frac{\partial \varphi}{\partial x} \cdot
      \left.\frac{\partial X}{\partial s}\right|_{s = 0}
      + \frac{\partial \varphi}{\partial u} \cdot
      \left.\frac{\partial U}{\partial s}\right|_{s = 0}
      + \varphi \frac{d}{dt} \left.\frac{\partial T}{\partial s}\right|_{s =
      0} \, , \label{eq:ds0ii} \\
0 = \frac{\partial g}{\partial t}
      \left.\frac{\partial T}{\partial s}\right|_{s = 0}
      + \frac{\partial g}{\partial x} \cdot
      \left.\frac{\partial X}{\partial s}\right|_{s = 0}
      + \frac{\partial g}{\partial u} \cdot
      \left.\frac{\partial U}{\partial s}\right|_{s = 0}
      + g \frac{d}{dt} \left.\frac{\partial T}{\partial s}\right|_{s =
      0} \, . \label{eq:ds0iii}
\end{gather}
Multiplying \eqref{eq:ds0i} by $\psi_0$, \eqref{eq:ds0ii} by $\psi(t)$,
and \eqref{eq:ds0iii} by $\lambda$, we can write:
\begin{multline}
\label{eq:joined}
\psi_0 \left(\frac{\partial L}{\partial t}
      \left.\frac{\partial T}{\partial s}\right|_{s = 0}
      + \frac{\partial L}{\partial x} \cdot
      \left.\frac{\partial X}{\partial s}\right|_{s = 0}
      + \frac{\partial L}{\partial u} \cdot
      \left.\frac{\partial U}{\partial s}\right|_{s = 0}
      + L \frac{d}{dt} \left.\frac{\partial T}{\partial s}\right|_{s =
      0}\right) \\
+ \psi(t) \cdot \left(\frac{\partial \varphi}{\partial t}
      \left.\frac{\partial T}{\partial s}\right|_{s = 0}
      + \frac{\partial \varphi}{\partial x} \cdot
      \left.\frac{\partial X}{\partial s}\right|_{s = 0}
      + \frac{\partial \varphi}{\partial u} \cdot
      \left.\frac{\partial U}{\partial s}\right|_{s = 0}
      + \varphi \frac{d}{dt} \left.\frac{\partial T}{\partial s}\right|_{s =
      0} - \frac{d}{dt} \left.\frac{\partial X}{\partial s}\right|_{s = 0}\right)\\
+ \lambda \cdot \left(\frac{\partial g}{\partial t}
      \left.\frac{\partial T}{\partial s}\right|_{s = 0}
      + \frac{\partial g}{\partial x} \cdot
      \left.\frac{\partial X}{\partial s}\right|_{s = 0}
      + \frac{\partial g}{\partial u} \cdot
      \left.\frac{\partial U}{\partial s}\right|_{s = 0}
      + g \frac{d}{dt} \left.\frac{\partial T}{\partial s}\right|_{s =
      0} \right) = 0 \, .
\end{multline}
According to the maximality condition
of the Pontryagin maximum principle, the function
\begin{multline*}
\psi_0 L\left(t,x(t),U\left(t,x(t),u(t),s\right)\right)
+ \psi(t) \cdot \varphi\left(t,x(t),U\left(t,x(t),u(t),s\right)\right)\\
+ \lambda \cdot g\left(t,x(t),U\left(t,x(t),u(t),s\right)\right)
\end{multline*}
attains an extremum for $s = 0$. Therefore
\begin{equation*}
\psi_0 \frac{\partial L}{\partial u} \cdot
\left.\frac{\partial U}{\partial s}\right|_{s = 0}
+ \psi(t) \cdot \frac{\partial \varphi}{\partial u} \cdot
\left.\frac{\partial U}{\partial s}\right|_{s = 0}
+ \lambda \cdot \frac{\partial g}{\partial u} \cdot
\left.\frac{\partial U}{\partial s}\right|_{s = 0} = 0
\end{equation*}
and \eqref{eq:joined} simplifies to
\begin{multline}
\label{eq:applASandDHdt}
\psi_0 \left(\frac{\partial L}{\partial t}
      \left.\frac{\partial T}{\partial s}\right|_{s = 0}
      + \frac{\partial L}{\partial x} \cdot
      \left.\frac{\partial X}{\partial s}\right|_{s = 0}
      + L \frac{d}{dt} \left.\frac{\partial T}{\partial s}\right|_{s =
      0}\right) \\
+ \psi(t) \cdot \left(\frac{\partial \varphi}{\partial t}
      \left.\frac{\partial T}{\partial s}\right|_{s = 0}
      + \frac{\partial \varphi}{\partial x} \cdot
      \left.\frac{\partial X}{\partial s}\right|_{s = 0}
      + \varphi \frac{d}{dt} \left.\frac{\partial T}{\partial s}\right|_{s =
      0} - \frac{d}{dt} \left.\frac{\partial X}{\partial s}\right|_{s = 0}\right)\\
+ \lambda \cdot \left(\frac{\partial g}{\partial t}
      \left.\frac{\partial T}{\partial s}\right|_{s = 0}
      + \frac{\partial g}{\partial x} \cdot
      \left.\frac{\partial X}{\partial s}\right|_{s = 0}
      + g \frac{d}{dt} \left.\frac{\partial T}{\partial s}\right|_{s =
      0} \right) = 0 \, .
\end{multline}
From the adjoint system
$\dot{\psi} = -\frac{\partial H}{\partial x}$ and the equality
\eqref{eq:dHdt}, we know that
\begin{gather*}
\dot{\psi} = - \psi_0 \frac{\partial L}{\partial x}
- \psi \cdot \frac{\partial \varphi}{\partial x}
- \lambda \cdot \frac{\partial g}{\partial x} \\
\frac{d}{dt} H = \psi_0 \frac{\partial L}{\partial t}
+ \psi \cdot \frac{\partial \varphi}{\partial t}
+ \lambda \cdot \frac{\partial g}{\partial t} \, ,
\end{gather*}
and one concludes that \eqref{eq:applASandDHdt} is equivalent to
\begin{equation*}
\frac{d}{dt} \left(\psi(t) \cdot
\left.\frac{\partial X}{\partial s}\right|_{s = 0}
- H \left.\frac{\partial T}{\partial s}\right|_{s = 0}\right)
= 0 \, .
\end{equation*}
The proof is complete.
\end{proof}

We now introduce the notion of unimprovable or Pareto conservation law.
\begin{definition}
An equation $C\left(t,x(t),u(t),\psi(t),\lambda\right) = \text{constant}$,
valid in $t \in [a,b]$ for any quadruple
$\left(x(\cdot),u(\cdot),\psi(\cdot),\lambda\right)$
satisfying the ``general theorem of optimal control''
(Theorem~\ref{th:GenThOptCont}),
is called an \emph{unimprovable conservation law}
or a \emph{Pareto conservation law} for problem $(P)$.
\end{definition}
Given the relation between problems $(P_1)$ and $(P)$
(\textrm{cf.} Section~\ref{sec:VecValOCP}), we obtain from
Theorem~\ref{r:MC:consLaw:IC} the following corollary.
\begin{theorem}[Noether theorem for vector-valued optimal control systems]
\label{r:MC:consLaw:PO}
If there exists a $C^2$-smooth one-parameter group of transformations
\begin{gather*}
h^s : [a,b] \times \mathbb{R}^n \times \mathbb{R}^r \rightarrow
       \mathbb{R} \times \mathbb{R}^n \times \mathbb{R}^r \, , \\
h^s(t,x,u) = \left(T(t,x,u,s), X(t,x,u,s), U(t,x,u,s)\right) \, , \\
s \in (-\varepsilon, \varepsilon)  \, , \, \varepsilon > 0 \, ,
\end{gather*}
with $h^0(t,x,u) = (t,x,u)$ for all
$(t,x,u) \in [a,b] \times \mathbb{R}^n \times \mathbb{R}^r$,
and satisfying
\begin{gather}
\frac{d}{dt} X\left(t,x(t),u(t),s\right)
= \varphi \circ h^s\left(t,x(t),u(t)\right)
\frac{d}{dt} T\left(t,x(t),u(t),s\right) \, ,  \label{eq:INV:PVV:X}\\
L\left(t,x(t),u(t)\right)
= L \circ h^s\left(t,x(t),u(t)\right)
\frac{d}{dt} T\left(t,x(t),u(t),s\right)  \, , \label{eq:INV:PVV:L}
\end{gather}
($L = \left(L_1,\ldots,L_N\right)$) then,
\begin{multline}
\label{eq:MC:consLaw:PO}
\psi(t) \cdot \frac{\partial}{\partial s}
\left.X\left(t,x(t),u(t),s\right)\right|_{s = 0} \\
- \mathcal{H}(t,x(t),u(t),\psi(t),\lambda) \frac{\partial}{\partial s}
\left.T\left(t,x(t),u(t),s\right)\right|_{s = 0} = \text{const}
\end{multline}
is an unimprovable conservation law for problem $(P)$,
with $\mathcal{H}$ the Hamiltonian \eqref{eq:Hamiltonian:P}
associated to the problem $(P)$.
\end{theorem}

\begin{remark}
Theorems~\ref{r:MC:consLaw:IC} and \ref{r:MC:consLaw:PO} are
still valid in the situation where the boundary values
of the state variables and/or the initial-terminal
instants of time ($a$, $b$) are not fixed.
We have considered conditions \eqref{eq:bc} and fixed
both initial and terminal times, only to simplify
the presentation of the Pontryagin maximum principle:
initial and terminal transversality conditions are not relevant
in the proof of our Noether-type theorems.
\end{remark}

In the next section we illustrate Theorem~\ref{r:MC:consLaw:PO}
with an example of five state variables ($n=5$), two controls ($r=2$),
and two functionals to minimize ($N = 2$).


\section{Example for the Flight of a Pilotless Aircraft}

We borrow from \cite[\S 3.4]{MR81k:93001} the problem
of optimizing a vector functional with two components,
representing fuel expenditure ($I_1$) and flight-time ($I_2$),
\begin{equation*}
I_1 = \int_0^T u_1(t) \mathrm{d}t \, , \quad
I_2 = \int_0^T 1 \mathrm{d}t \, ,\\
\end{equation*}
subject to a dynamical control system representing
the motion of a pilotless aircraft
\begin{equation*}
\begin{cases}
\dot{x}_1(t) &= x_3(t) \, ,\\
\dot{x}_2(t) &= x_4(t)  \, ,\\
\dot{x}_3(t) &= c_1 \frac{u_1(t)}{x_5(t)} \cos(u_2(t))  \, ,\\
\dot{x}_4(t) &= c_1 \frac{u_1(t)}{x_5(t)} \sin(u_2(t)) - c_2  \, ,\\
\dot{x}_5(t) &= - u_1(t) \, .
\end{cases}
\end{equation*}
Here $x_1$ is the range of the aircraft;
$x_2$ the altitude;
$x_3$ the horizontal component of the velocity;
$x_4$ the vertical component of the velocity;
$x_5$ the mass of the aircraft (which depends of its fuel);
$u_1$ the rate of fuel consumption;
$u_2$ the thrust angle relative to the horizontal;
$c_1$ and $c_2$ given constants. A full description
of the model, and a complete analysis of its solution,
is found on \cite[\S 3.4]{MR81k:93001}. Our objective
here is to obtain a non-trivial unimprovable conservation law
for the problem, with the help of Theorem~\ref{r:MC:consLaw:PO}.
About the model, it is enough for our purposes to say
that there are physical constraints on the control values,
under which makes sense to consider $\tan(u_2)$
(\textrm{cf.} \cite[(3.42)]{MR81k:93001}).
Two trivial unimprovable conservation laws are $\psi_1(t) = \text{const}$
(obtained from Theorem~\ref{r:MC:consLaw:PO} choosing
$T = t$, $X_1 = x_1 + s$, $X_i = x_i$, $i = 2,\ldots,5$, $U_j = u_j$, $j = 1,2$),
and $\psi_2(t) = \text{const}$ (obtained from Theorem~\ref{r:MC:consLaw:PO} choosing
$T = t$, $X_2 = x_2 + s$, $X_i = x_i$, $i = 1,3,4,5$, $U_j = u_j$, $j = 1,2$).
We claim that
\begin{equation}
\label{eq:CL:FPA}
\psi_1 x_1(t) + 2 \psi_2 x_2(t) + \psi_3(t) x_3(t) + 2 \psi_4(t) x_4(t) = \text{const}
\end{equation}
is also an unimprovable conservation law for the problem. We
remark that \eqref{eq:CL:FPA} is non-trivial, and difficult to
obtain without Theorem~\ref{r:MC:consLaw:PO}. To prove it with the
help of Theorem~\ref{r:MC:consLaw:PO}, one just need to show that
the problem is invariant (satisfies conditions
\eqref{eq:INV:PVV:X} and \eqref{eq:INV:PVV:L}) with $T = t$, $X_1
= \mathrm{e}^s x_1$, $X_2 = \mathrm{e}^{2s} x_2$, $X_3 =
\mathrm{e}^s x_3$, $X_4 = \mathrm{e}^{2s} x_4$, $X_5 = x_5$, $U_1
= u_1$, and $U_2 = \arctan\left(\mathrm{e}^s \tan u_2\right)$
($\sin U_2 = \mathrm{e}^{2s} \sin u_2$, $\cos U_2 = \mathrm{e}^{s}
\cos u_2$). This is done by direct calculations ($\frac{d}{dt} T = 1$):
\begin{equation*}
\begin{split}
\frac{d}{dt} X_1 &= \mathrm{e}^{s} \dot{x}_1 = \mathrm{e}^{s} x_3 = X_3 \frac{d}{dt} T \, , \\
\frac{d}{dt} X_2 &= \mathrm{e}^{2s} \dot{x}_2 = \mathrm{e}^{2s} x_4 = X_4 \frac{d}{dt} T \, , \\
\frac{d}{dt} X_3 &= \mathrm{e}^{s} \dot{x}_3 = c_1 \frac{u_1}{x_5}\mathrm{e}^{s} \cos u_2
                                             = c_1 \frac{U_1}{X_5} \cos U_2 \frac{d}{dt} T \, , \\
\frac{d}{dt} X_4 &= \mathrm{e}^{2s} \dot{x}_4 = c_1 \frac{u_1}{x_5}\mathrm{e}^{2s} \sin u_2 - c_2
                  = \left(c_1 \frac{U_1}{X_5} \sin U_2 - c_2\right) \frac{d}{dt} T \, , \\
\frac{d}{dt} X_5 &=  \dot{x}_5 = - u_1 = - U_1 \frac{d}{dt} T \, ,
\end{split}
\end{equation*}
and equations \eqref{eq:INV:PVV:X} are verified;
\begin{equation*}
\begin{split}
L_1 &= u_1 = U_1 \frac{d}{dt} T \, , \\
L_2 &= 1 = \frac{d}{dt} T \, ,
\end{split}
\end{equation*}
and equations \eqref{eq:INV:PVV:L} are also satisfied.
Equality \eqref{eq:MC:consLaw:PO} takes then form \eqref{eq:CL:FPA}.


\section*{Acknowledgements}

The author acknowledges the Control Systems Department
of the Georgian Technical University in Tbilisi,
for the invitation to visit Georgia on September 2004;
for the reference \cite{MR81k:93001}; and for giving
him the opportunity to learn about
the excellent and interesting work
which is done in Tbilisi on the problems of optimal control,
with the help of symmetry and conservation laws,
and the applications in concrete fields
of seismology, energetic chemistry, and metallurgy.
The author is particularly grateful to Valida Sesadze
and Tamuna Kekenadze.



\end{document}